\documentstyle[11pt]{article}

\title{\Large \bf  A note on Closure Operators in Category of Groups}

\author{Vishvajit V S Gautam\thanks{E-mail : vishvajit@imsc.res.in; gautamvvs@yahoo.com }\\
\small \sl  The Institute of Mathematical Sciences\\
\small \sl  CIT Campus Taramani, Chennai - 600113 INDIA}
\date{ }
\begin{document}
\newtheorem{defi}{\bf Definition}[section]
\newtheorem{th}[defi]{\bf Theorem}
\newtheorem{pro}[defi]{\bf Proposition}
\newtheorem{lemma}[defi]{\bf Lemma}
\newtheorem{coro}[defi]{\bf Corollary}
\newtheorem{ram}[defi]{\bf Remark}
\newtheorem{pn}[subsubsection]{\bf Proposition}
\newtheorem{zm}[subsubsection]{\bf Theorem}
\newtheorem{cy}[subsubsection]{\bf Corollary}
\newtheorem{df}[subsubsection]{\bf Definition}
\newtheorem{rk}[subsubsection]{\bf Remark}
\maketitle
\begin{abstract}
We give some  applications of closure operators in category 
of groups and link them with the {\em join problem} of subnormal subgroups.\\


\end{abstract}

Notion of  closure operators(operations,  systems,  functions, 
relations) is known to us from algebra. logic, lattice theory and topology.
Categorical view of closure operators play an important role in various
branches of mathematics. In an arbitrary category $\cal X$ with suitable
axiomatically defined notion of subobjects, a (categorical) 
closure operator $c$ is defined to be a family $(c_{X})_{X \in {\cal X}}$ 
satisfying the properties of extension, monotonicity and continuity.
 Closure operators are proved to be
useful in study of Galois equivalence between certain factorization systems. 
In category of {\bf R-mod} of R-modules closure operators correspond to 
preradicals. For more details, see [2]. 
In this article we establish a few results and  examples in
the category of groups  by means of closure operators. In section 2, 
Theorem 2.1.3 provides an interesting link between the {\em join} problem
of two subnormal subgroups and additive closure operators defined on 
{\bf Grp} the category of groups. In subsection 2.2 we use the notion of
 closure operator induced by a subcategory $\cal A$ of {\bf Ab} the
 category of abelian groups to characterize the homomorphisms between the
quotient group $G/H$ and  a group  $A \in \cal A$. This characterization
provide  useful methods for determining the relation 
between epimorphisms 
and surjective homomorphisms in many algebraic categories.

Subcategories are always assumed to be full and isomorphism closed.

\section{Preliminaries}
Throughout this paper we consider a category $\cal X$ and a fixed class $\cal M$ 
of monomorphisms in $\cal X$ which contains all isomorphisms of $\cal X$. It is 
assumed that\\
$\bullet$ $\cal M$ is closed under composition;\\
$\bullet$ $\cal X$ is finite $\cal M$-complete.\\

A {\em closure operator} $c$ on the category $\cal X$ with respect to class 
$\cal M$ of subobjects is given by a family $c = (c_{X})_{X \in \cal X}$ of 
maps $c_{X} : {\cal M}/X \longrightarrow {\cal M}/X$ such that for every 
$X \in \cal X$ \\
1.  $m  \leq  c(m)$;  2.   $m \leq m' \Rightarrow  c (m) \leq c (m')$;
and  3. for  every $f:X \longrightarrow Y$ and  $m \in {\cal M}/X$,
\hspace{.5in} $ f(c_{X}(m)) \leq c_{Y}(f(m))$.\\
\\
For each $m \in {\cal M}$  we denote by $c(m)$ the c-closure of $m$.\\
An $\cal M$-morphism   $m \in {\cal M}/X$ is  called $c$-closed
if  $m \cong c_{X}(m)$. A closure operator $c$ is said to be idempotent
if  $c(c(m)) \cong c(m)$. In case $ c(m \vee n) \cong c(m) \vee c(n)$ we say 
$c$ to be additive. An $\cal M$-subobject $m$ of $X$ is called $c$-dense in
$X$ if $c_{X}(m) \cong 1_{X}$.
\\
For a subcategory $\cal A$ of $\cal X$, a morphism $f:X \longrightarrow Y$
is an $\cal A$-regular monomorphism  if it is the equalizer of two morphisms
 $h,k:Y \longrightarrow A$  with  $A \in {\cal A}$.

 Let   $\cal M$  contain  the  class  of  regular monomorphisms of $\cal X$.
 For $m:M \longrightarrow X$ in $\cal M$ define 
 \[c_{\cal A}(m) = \bigwedge \{ r \in {\cal M} \mid r \geq m \,\,\,\, and\,\,
 r \,\,\,is \,\, {\cal A}-regular \} \]
which  is  a closure  operator  of $ {\cal X}$.  These closure operators
are called regular  and $c_{\cal A}(m)$  is  called  the  $\cal A$-closure
of $m$. In case ${\cal A} = {\cal X}$ we denote $c_{\cal A}(m)$ by $c(m)$.

\section{Closure operators in category of groups}
In this section we will see application of closure operators in category of
groups.

\subsection{ }
Let ${\cal X} = $ {\bf Grp} the category of groups and 
let $\cal M$ be the class of all monomorphisms of $\cal X$. In this case clearly
$\cal X$ is finite $\cal M$-complete. For an object $G$ of {\bf Grp}, 
${\cal M}/G$ can be identified with the set of all subgroups of $G$.\\
Let $H$ be a subgroup of $G$, we define
\[c_{G}(H) = <g^{-1}hg \mid h \in H, g \in G>\]
the least normal subgroup of $G$ containing $H$.\\
It is easy to prove that $c_{G}$ is a closure operator on $\cal X$ which 
is also an idempotent operator.\\
Let $[G, G] = <xyx^{-1}y^{-1} \mid x, \, y \in G>$ be the commutator subgroup 
of $G$. Define
\[ c'_{G}(H) = [G, G] \cdot H \]
This gives a closure operator on {\bf Grp} which is normal in $G$. For trivial
subgroup $(e)$ of $G$, $c'_{G}(e) = [G, G]$ while $c_{G}(e) = (e)$.\\
A {\em preradical} {\bf r} in {\bf Grp} is the subfunctor of the identity
functor in {\bf Grp}. For $G \in$ {\bf Grp}, {\bf r}($G$) is a normal
subgroup of $G$. We define two more closure operators on {\bf Grp} as follows:
$$c''_{G}(H) = H \cdot {\bf r}(G)  =  H \cdot {\bf r}(G) $$
and
$${c'''_{G}(H) = \pi^{-1}} ({\bf r}(G)/{c_{G}(H) })$$
where $\pi : G \longrightarrow {G/{c_{G}(H)}}$ is the canonical projection.\\
Closure operator $c''_{G}$ (in general) is not normal in $G$, but $c'''_{G}$
is always normal in $G$. (see [2])\\
We can observe that closure of a subgroup $H$ of $G$ can be converted to a
normal closure of $H$ in $G$ and vise-versa. For example,\\
\[c_{G}(c''_{G}(H)),\,\,\, c'_{G}(c''_{G}(H)),\,\,\, c'''_{G}(c''_{G}(H)) \] are
normal in $G$, but \[c''_{G}(c_{G}(H)),\,\,\, c''_{G}(c'_{G}(H)),\,\,\, 
c''_{G}(c'''_{G}(H))\] are ( in general) not normal in $G$.\\
\\
Next result is obvious.
\begin{pn} 
Let $G,H$ be objects in ${\cal X} = $ {\bf Grp}. Let $f: G \longrightarrow H$ 
be a non-zero homomorphism from $G$ to a simple group $H$.  $f$ 
is onto if and only if $c_{f(G)}(K) = f(G)$ for all subgroups $(e) \neq K$ of 
$f(G)$.
\end{pn}

Recall that a subgroup $H$ of a group $G$ is said to be {\em subnormal} in $G$
if there are a non-negative integer $m$ and a series
\[ H = H_{m} \triangleleft H_{m-1} \triangleleft H_{m-2} \triangleleft
\ldots  \triangleleft H_{0} = G \]
of subgroups of $G$. In this situation we write $H \,\,sm\,\, G$ and 
$H {\triangleleft^{m}} G$. The smallest such $m$ is called the defect of
subnormal subgroup $H$ of $G$.\\
\\
In finite group theory subnormal subgroups are precisely those subgroups
which occur as terms of composition series, the factors of which are of paramount
importance in describing a group's structure. In 1939 H. Wielandt proved
the celebrated {\em join} theorem for finite groups. Twenty year latter
H. Zassenhaus showed that Wielandt's {\em join} theorem can fail to hold in
infinite groups. The determination of interesting necessary and sufficient 
conditions for a {\em join} of subnormal subgroups (i.e. the subgroup 
generated by two subnormal subgroups) to be subnormal is probably the 
most important unsolved problem
in this area of group theory [4]. In following we establish an interesting link
between subnormal subgroups and closure operators which provides a solution of
subnormal subgroups {\em join} problem.\\
\\
Let $c_{G}(H)$ denote the  closure of $H$ in $G$. Set $H_{0} = G$,
$H_{1} = c_{H_{0}}(H)$ , $H_{2} = c_{H_{1}}(H)$, $\ldots$, 
$H_{m+1} = c_{H_{m}}(H)$.

\begin{pn}
([4]) Let $H$ be a subgroup of $G$. Then
$H {\triangleleft^{m}} G$ if and only if $H$ coincides with its mth normal
closure in $G$.
\end{pn}

Denote by $T = <H,K>$ group generated by two subnormal subgroups and
$T_{m,n} = <H_{m},K_{n}>$ for $m,n = 0,1,2, \ldots $ where
$H_{m+1} = c_{H_{m}}(H)$ and $K_{n+1} = c_{K_{n}}(K)$.

\begin{zm}
Let $H$ and $K$ be subnormal subgroups of a group $G$
in $\cal X$. Then following implications hold.
\begin{enumerate}
\item The class ${\cal M}{^{c}}{_{G}}$ of $c$-closed elements in ${\cal M}/G$
is closed under binary suprema for every object $G$ in $\cal X$.
\item $c$ is additive.
\item $c_{T_{m,n}}(T)$ is subnormal in $G$.
\end{enumerate}
1. $\Longleftrightarrow$ 2. $\Longrightarrow$ 3. 
\end{zm}
{\bf Proof.} (sketch) 1 $\Longleftrightarrow$ 2  straightforward.\\
2 $\Longrightarrow$ 3 . Let
\[ H \leq \ldots c_{H_{m}}(H) \triangleleft c_{H_{m-1}}(H) \triangleleft 
c_{H_{m-2}}(H) \triangleleft \ldots  \triangleleft c_{H_{1}}(H) 
\triangleleft c_{G}(H) \triangleleft G \]
and
\[ K \leq \ldots c_{K_{n}}(K) \triangleleft c_{K_{n-1}}(K) \triangleleft 
c_{K_{n-2}}(K) \triangleleft \ldots  \triangleleft c_{K_{1}}(K) 
\triangleleft c_{G}(K) \triangleleft G \]
be the series of $H$ and $K$ respectively.\\
Since $c$ is additive, therefore we have\\
$c_{G} (H \vee K) = c_{G}(H) \vee c_{G}(K)$ i.e. 
$c_{G}(<H,K>) = <c_{G}(H), c_{G}(K)>$.\\
Also we have\\
\[c_{G}(<c_{H_{i-1}}(H), c_{K_{j-1}}(K)>) = 
<c_{G}(c_{H_{i-1}}(H)), c_{G}(c_{K_{j-1}}(K))>\]
for all $i,\,\, 1 \leq i \leq m$ and $j,\,\, 1 \leq i \leq n$, and for all
group $G$ in $\cal X$.\\
Clearly we have
\[ c_{T_{m, n}}(<H, K>) 
\subseteq c_{T_{m, n}}(<c_{H_{m}}(H), c_{K_{n}}(K)>) \] 
\[ \subseteq c_{T_{m-1, n-1}}(<c_{H_{m-1}}(H), c_{K_{n-1}}(K)>)\]
\[ \subseteq c_{T_{m-2, n-2}}(<c_{H_{m-2}}(H), c_{K_{n-2}}(K)>)\] 
\[ \subseteq \ldots  \subseteq c_{T_{1, 1}}(<c_{H_{1}}(H), c_{K_{1}}(K)>) 
\subseteq G. \]

By the additivity of $c$ and the above expression we get
\[ c_{T_{m, n}}(T) \triangleleft <c_{H_{m-1}}(H), c_{K_{n-1}}(K)> \triangleleft <c_{H_{m-2}}(H), c_{K_{n-2}}(K)> \]
\[ \triangleleft \ldots \triangleleft <c_{H_{1}}(H), c_{K_{1}}(K)>
\triangleleft \,\, G.\]
This proves that $ c_{T_{m, n}}(T)$ is subnormal in $G$. $\Box$

\begin{rk}
 If $c$ is a normal closure operator in $G$, then it
satisfy the conditions 1 and 2 of the Theorem 2.1.3. This means that
{\em join} problem of two subnormal subgroups is reduced to find the suitable
conditions when $H = c_{H_{m}}(H)$ and $K = c_{K_{n}}(K)$ (cf. Prop. 2.1.2). 
\end{rk}

\begin{rk}
Normal subgroups are not only stable under intersection, but also under
arbitrary join in the subgroup lattice. Therefore all the closure operators 
which are normal as a subgroup of $G$ are additive. 
We observe that the normal subgroups  produced by $c'_{G}$, $c''_{G}$ and $c'''_{G}$
will be larger than of the normal subrgoups produced by $c_{G}$.
So there will be possibility of a fast  termination of subnormal series induced by these operators.
\end{rk}
As a corollary to above theorem we have following result of Wielandt
(cf. [4])
\begin{cy}
 If $H$ and $K$ are subnormal subgroups of of a
finite group $G$ in $\cal X$. Then $T = <H,K>$ is subnormal in G.
\end{cy}

\subsection{ }
Let $\cal A$ be a subcategory of ${\cal X} = $ {\bf Ab} the
category of abelian groups and let $\cal M$ be the class of all monomorphisms
of $\cal X$. Notice that in {\bf Ab} and {\bf Grp} the strong monomorphisms
coincide with monomorphisms. In this case $\cal X$ is $\cal M$-complete.\\
For an object $G \in${\bf Ab}, ${\cal M}/G$ can be identified with the set 
of all subgroups of G. 

\begin{zm}
 Let $H$ be a subgroup of $G$. \\
{\bf (a)} (cf. [1]) $H$ is $\cal A$-dense in $G$, i.e., $c_{\cal A}(H) = G$
if and only if $Hom(G/H , A) = (0)$ for every $A \in \cal A$.\\
{\bf (b)} $H$ is $\cal A$-closed, i.e., $c_{\cal A}(H) = H$ if and only if
$Hom(G'/H , A) \not= (0)$ for some $A \in \cal A$ and for every
non-zero subgroup $G'/H$ of $G/H$.\\
{\bf (c)} Let $f:G \longrightarrow T$ be an $\cal A$-morphism. The morphism 
$f$ is not epic if and only if $c_{\cal A}(f(G)) = f(G)$.\\
In particular $c_{\cal A}(f(G)) = f(G)$ implies f is not surjective.
\end{zm}
{\bf Proof.} (a) (cf. [1])\\
(b) Without loss of generality we assume that $H$ is not $\cal A$-dense in $G$.
Let $G'/H$ be a non-zero subgroup of $G/H$. Since $H$ is not $\cal A$-dense 
in $G$, we have $Hom(G/H, A) \not= (0)$ for some $A \in \cal A$. Let 
$0 \not= f \in Hom(G/H, A)$. One can get a non-zero morphism
$f':G'/H \longrightarrow A$ in obvious sense i.e., $f' = f \cdot j$ where
$j:G'/H \longrightarrow G/H$ is just the inclusion map, which implies that
$Hom(G'/H, A) \not= (0)$.\\
Conversely, suppose $H$ a proper subgroup of $c_{\cal A}(H)$. Since closure
of  $H$ is $c_{\cal A}(H)$, $H$ is dense in $c_{\cal A}(H)$ (treating $H$
as a subobject fo $c_{\cal A}(H)$). This implies  
$Hom({c_{\cal A}(H)}/H, A) = (0)$,
but this contradicts our hypothesis, therefore we must have 
$c_{\cal A}(H) = H$.\\
\\
(c) Since $f(G)$ is closed in $T$ implies $Hom(T/f(G), A) \not= (0)$ for
some $A \in \cal A$. This implies $f$ is not an epimorphism.\\
Conversely, if $f$ is not an epimorphism implies $c_{\cal A}(f(G)) \not= T$
(cf. [1]) which gives $Hom(T'/f(G), A) \not= (0)$ for some $A \in \cal A$
and for every non-zero subgroup $T'/f(G)$ of $T/f(G)$. But then from (b)
we get $c_{\cal A}(f(G)) = f(G)$. $\Box$

\begin{rk}
 Above results can be used in case of following 
subcategories of {\bf Ab}. (cf. [3])\\
1. Category of torsion free abelian groups; 2. category of reduces abelian
groups; 3. category of free abelian groups; 4. category of  topological
abelian groups, etc.
\end{rk}
\vspace{.5in}
{\Large References}
\newline [1]   G. Castellini, {\em Closure operators, monomorphisms and 
epimorphisms
in categories of groups}, Cahiers de Topologie et Geometric Differentielie
Categoriques Vol. XXVII -2(1986)151-167.
\newline [2]   D.  Dikranjan  and  W.  Tholen, {\em Categorical  structure of 
closure  operators  with  applications  to  Topology, Algebra and 
Discrete Mathematics} (Kluwer, Dordrecht, 1994).
\newline [3]   L. Fuchs, {\em Infinite abelian groups Vol. 1} (Academic Press, 
New York, 1970)
\newline [4]   J.C. Lennox and S.E. Stonehewer, {\em Subnormal subgroups 
of groups} (Clarendon Press, Oxford, 1987). 

\end{document}